\RequirePackage{ifpdf}
\ifpdf 
\documentclass[pdftex]{sigma}
\else
\documentclass{sigma}
\fi

\numberwithin{equation}{section}

\newcommand{\x}{\times}
\newcommand{\D}{\Delta}
\newcommand{\G}{\Gamma}
\newcommand{\p}{\partial}

\newcommand{\bx}{{\mathbf x}}
\newcommand{\bm}{{\mathbf m}}
\newcommand{\bu}{{\mathbf u}}

\begin{document}

\allowdisplaybreaks

\renewcommand{\PaperNumber}{090}

\FirstPageHeading

\ShortArticleName{On a Family of 2-Variable Orthogonal Krawtchouk Polynomials}

\ArticleName{On a Family of 2-Variable Orthogonal Krawtchouk \\ Polynomials}

\Author{F. Alberto GR\"UNBAUM~$^\dag$ and Mizan RAHMAN~$^\ddag$}

\AuthorNameForHeading{F.A.~Gr\"unbaum and M.~Rahman}

\Address{$^\dag$~Department of Mathematics, University of California, Berkeley, CA 94720, USA}
\EmailD{\href{mailto:grunbaum@math.berkeley.edu}{grunbaum@math.berkeley.edu}}
\URLaddressD{\url{http://www.math.berkeley.edu/~grunbaum/}}

\Address{$^\ddag$~Department of Mathematics and Statistics, Carleton University, Ottawa, Canada, K1S~5B6}
\EmailD{\href{mailto:mrahman@math.carleton.ca}{mrahman@math.carleton.ca}}

\ArticleDates{Received July 25, 2010, in f\/inal form December 01, 2010;  Published online December 07, 2010}

\Abstract{We give a hypergeometric proof involving a family of 2-variable
Krawtchouk polynomials that were obtained earlier by Hoare and Rahman [{\it SIGMA} {\bf 4} (2008), 089, 18~pages]
as a limit of the $9-j$ symbols of quantum angular momentum theory, and
shown to be eigenfunctions of the transition probability kernel
corresponding to a ``poker dice'' type probability model. The proof in this
paper derives and makes use of the necessary and suf\/f\/icient conditions of
orthogonality in establishing orthogonality as well as indicating their
geometrical signif\/icance. We also derive a 5-term recurrence relation
satisf\/ied by these polynomials.}


\Classification{33C45}

\section{Introduction}
\label{sec1}

It was in the SIDE8 meeting in St-Adele near Montreal that one of us (MR) presented a paper reporting the discovery, by his co-author Michael Hoare and himself [2008], of a ``new'' system of 2-variable Krawtchouk polynomials, orthogonal with respect to a trinomial distribution.  The motivation of their paper was to f\/ind eigenvalues and eigenfunctions of the transition probability kernel:
\begin{gather}
K_A(j_1,j_2;i_1,i_2)  = \sum_{k_1=0}^{\min(i_1,j_1)} \sum_{k_2=0}^{\min(i_2,j_2)} b(k_1,i_1;\alpha_1)b(k_2,i_2;\alpha_2) \nonumber\\
\phantom{K_A(j_1,j_2;i_1,i_2)=}{} \x b_2(j_1-k_1,j_2-k_2;N-k_1-k_2;\beta_1,\beta_2),\label{eq1.1}
\end{gather}
where $b(x,N;p) = \binom{N}{x} p^x(1-p)^{N-x}$ is the binomial distribution, while
\begin{gather*}
b_2(x,y;N;p,q) = \binom{N}{x,y} p^xq^y(1-p-q)^{N-x-y},
\end{gather*}
is the trinomial, both normalized to $1$.  In \eqref{eq1.1} the parameters $(\alpha_1,\alpha_2,\beta_1,\beta_2)$ are probabilities of a two-step cumulative Bernoulli process, and hence necessarily in $(0,1)$, while $(i_1,i_2)$ and $(j_1,j_2)$ represent the initial and f\/inal states of the process.  Borrowing a result from the angular momentum theory of quantum mechanics the authors of \cite{HR4} were able to show that the $2$-dimensional Krawtchouk polynomials
\begin{gather}
\label{eq1.3}
{\underset{0 \le i+j+k+l \le N}{\sum_i\sum_j\sum_k\sum_l}} \frac {(-m)_{i+j}(-n)_{k+l}(-x)_{i+k}(-y)_{j+l}}{i!j!k!l!(-N)_{i+j+k+l}} u_1^iv_1^ju_2^kv_2^l
\end{gather}
do indeed satisfy the requirements for them being the eigenfunctions of \eqref{eq1.1}, where $(x,y)$ represents the state-variable and $(m,n)$ the spectral parameters.  It goes without saying that for~\eqref{eq1.3} to be an orthogonal system with respect to a distribution of the form $b_2(x,y;N;\eta_1,\eta_2)$ the parameters $u_1$, $v_1$, $u_2$, $v_2$ must be related to $\eta_1$, $\eta_2$, as well as satisfy some additional conditions among themselves.  It was found in~\cite{HR4}, again with a cue from the physics literature that these conditions are all satisf\/ied provided the $u$'s and $v$'s are parametrized in the following way
\begin{alignat}{3}
& u_1 = \frac {(p_1+p_2)(p_1+p_3)}{p_1(p_1+p_2+p_3+p_4)}, \qquad && u_2 = \frac {(p_1+p_2)(p_2+p_4)}{p_2(p_1+p_2+p_3+p_4)}, & \nonumber\\
& v_1 = \frac {(p_1+p_3)(p_4+p_3)}{p_3(p_1+p_2+p_3+p_4)}, \quad && v_2 = \frac {(p_2+p_4)(p_3+p_4)}{p_4(p_1+p_2+p_3+p_4)}, & \label{eq1.4}
\end{alignat}
and consequently,
\begin{gather}
\label{eq1.5}
\eta_1 = \frac {p_1p_2(p_1+p_2+p_3+p_4)}{(p_1+p_2)(p_1+p_3)(p_2+p_4)}, \qquad
\eta_2 = \frac {p_3p_4(p_1+p_2+p_3+p_4)}{(p_2+p_4)(p_3+p_4)(p_1+p_3)}.
\end{gather}

For the origin of the work in \cite{HR4} the reader may consult \cite{CHR,HR1,HR2,HR3}.

Fortunately, in the audience, a very attentive listener, Masatoshi Noumi, was there to point out to (MR) that these polynomials are not only not new, but a special case of the multivariable generalization of the Gaussian hypergeometric function:
\begin{gather}
\label{eq1.6}
F_1^{(n)}(-\bx,-\bm;-N;\bu) = \sum \frac { \prod\limits_{i=1}^n (-x_i)_{\sum\limits_{j=1}^n \alpha_{ij}}  \prod\limits_{i=1}^n (-m_i)_{\sum\limits_{j=1}^n \alpha_{ji}}} {(-N)_{\sum\limits_{i,j} \alpha_{ij}}}  \frac {\prod u_{ij}^{\alpha_{ij}}}{\prod \alpha_{ij}!},
\end{gather}
where the $\alpha_{ij}$'s are nonnegative integers taking values from $0$ to $n$, such that $\sum\limits_{i,j}  \alpha_{ij} \le N$, $N$ being assumed as a nonnegative integer. Here we are following the notation in~\cite{M}.
In the original def\/inition of Aomoto and Gelfand \cite{AK,G}, $N$ need not be an integer, nor the~$x$'s and~$m$'s.  Furthermore, the space on which their functions are introduced was a bit more general, a generalization we do not need for our purposes.  To be sure, these authors' primary interest was not to look at~\eqref{eq1.6} as a multidimensional extension of the Krawtchouk polynomials:
\begin{gather}
\label{eq1.7}
P_n(x) = {}_2F_1\big(-x,-n;-N;\eta^{-1}\big),
\end{gather}
rather some structures that they contain.  H.~Mizukawa~\cite{M} proved that the functions in \eqref{eq1.6} are the zonal spherical functions on a very special class of Gelfand pairs made up of the ref\/lection groups $G(r,1,n)$ and the symmetric group $S_n$.
For a very nice account of a way to obtain many discrete orthogonal polynomials in terms of certain $(n+1,m+1)$ hypergeometric functions see also the work of Mizukawa and Tanaka~\cite{MT}. As mentioned above, we learned from Professor M.~Noumi that these functions give the multivariable Krawtchouk polynomials independently obtained in \cite{HR4}.  In a very recent paper by Mizukawa, see~\cite{M1}, he has established the ortho\-go\-nality of Krawtchouk polynomials in~$n$ variables by using very dif\/ferent techniques from the ones in this paper. See also the additional comment at the end of our paper.

The origin of the work in \cite{CHR,HR1,HR2,HR3,HR4} is the analysis of a very concrete probabilistic model, namely
``poker dice''. Its corresponding eigenfunctions are seen in~\cite{HR4} to be given in terms of a family of polynomials that, as indicated above, are now identif\/ied with the Gelfand--Aomoto polynomials. It is likely that this may be the f\/irst probabilistic application of the Gelfand--Aomoto polynomials. They may  also be applicable to other models in the
physical sciences.

One should point out that the hypergeometric functions involve both parameters as well as variables. Depending on the
issue at hand one can consider these functions as depending on one or the other set of ``variables''. This is already the
case in the classical one variable case as indicated by the expression~\eqref{eq1.7} above. The fact that these Krawtchouk polynomials (or the higher level Hahn polynomials) could be so useful in analyzing naturally appearing models in statistical mechanics was not anticipated in the classical book by W.~Feller~\cite{F}, where one can read about the Ehrenfest as well as the Bernoulli--Laplace models. For several applications of the Krawtchouk polynomials to several parts of mathematics see~\cite{Kra}. For a very good general guide to the f\/ield see~\cite{AAR}.

Our objective in this paper is less ambitious in one sense and more in another~-- namely, that we still restrict ourselves to the $n=2$ case, but not necessarily on the ref\/lection group, but to the general situations where the parameters $u_{ij}$'s will be determined by the requirement of orthogonality.  Hoare and Rahman~\cite{HR4} have done that problem, but we will approach it from a~dif\/ferent angle.  We will refrain from parametrizing the $u_{ij}$ right at the outset, instead looking for conditions they must satisfy among them in order that the $2$-variable polynomials:
\begin{gather}
F^{(2)}_1(-m_1,-m_2;-x_1,-x_2;-N;u_1,v_1,u_2,v_2) \nonumber\\
 \qquad := \sum \frac {(-m_1)_{i+j}(-m_2)_{k+l}(-x_1)_{i+k}(-x_2)_{j+l}}{i!j!k!l!(-N)_{i+j+k+l}} u_1^iv_1^ju_2^kv_2^l   \equiv P_{m_1,m_2}(x_1,x_2)\label{eq1.8}
\end{gather}
become orthogonal with respect to the trinomial:
\begin{gather*}
b_2(x_1,x_2;N;\eta_1,\eta_2) = \binom{N}{x_1,x_2} \eta_1^{x_1}\eta_2^{x_2}(1-\eta_1-\eta_2)^{N-x_1-x_2}.
\end{gather*}

It may be worth mentioning that a prior knowledge of this weight function is not essential since one could easily derive it by using the binomial generating function of the polynomials $P_{m,n}(x,y)$.

We now state the main results in the paper, namely \eqref{eq1.10} and \eqref{eq1.11} below.

In Sections~\ref{sec2}, \ref{sec3} and \ref{sec4} we will show that the necessary and suf\/f\/icient conditions of orthogonality are:
\begin{gather}
\text{(a)}  \quad \eta_1 u_1 + \eta_2 v_1 = 1, \nonumber\\
\text{(b)}  \quad \eta_1 u_2 + \eta_2 v_2 = 1, \label{eq1.10}\\
\text{(c)}  \quad \eta_1 u_1u_2 + \eta_2 v_1v_2 = 1,\nonumber
\end{gather}
with the $\eta$'s assumed to be given such that $0 < \eta_1,\eta_2 < 1$ and $\eta_1 + \eta_2 < 1$.  One can easily verify that these three conditions are all satisf\/ied by \eqref{eq1.4} and \eqref{eq1.5}.

However, one of the main reasons for going back to this problem is to f\/ind a $5$-term recurrence relation for \eqref{eq1.8}, which interestingly, is more easily found by using the $p$'s as in \eqref{eq1.4} and \eqref{eq1.5} than using \eqref{eq1.10} instead.  In Section \ref{sec5} we'll show that, if we denote \eqref{eq1.8} by $P_{m_1,m_2}(x_1,x_2)$ then
\begin{gather}
(N-m_1-m_2) \left\{ \frac {p_1p_3(p_2+p_4)(p_1+p_2+p_3+p_4)}{(p_1+p_3)(p_1p_4-p_2p_3)} (P_{m_1+1,m_2}(x_1,x_2) - P_{m_1,m_2}(x_1,x_2))\right. \nonumber\\
\qquad \left. {} - \frac {p_2p_4(p_1+p_3)(p_1+p_2+p_3+p_4)}{(p_2+p_4)(p_1p_4 - p_2p_3)} (P_{m_1,m_2+1}(x_1,x_2) - P_{m_1,m_2}(x_1,x_2))\right\} \nonumber\\
 \qquad {} + m_1 \frac {p_1p_4-p_2p_3}{p_1+p_3} (P_{m_1-1,m_2}(x_1,x_2) - P_{m_1,m_2}(x_1,x_2)) \nonumber\\
 \qquad {} - m_2 \frac {p_1p_4-p_2p_3}{p_2+p_4} (P_{m_1,m_2-1}(x_1,x_2) - P_{m_1,m_2}(x,y)) \nonumber\\
 \qquad {} = ((p_1+p_2)x_1 - (p_3+p_4)x_2) P_{m_1,m_2}(x_1,x_2).\label{eq1.11}
\end{gather}

This recursion relation is valid when the variables $x_1$, $x_2$ are nonnegative integers taking values whose sum is at most $N$.

If we insist on a dif\/ference operator in the variables $m_1$, $m_2$ with an eigenvalue that is linear in~$x_1$,~$x_2$ and involves only the four nearest neighbours of $m_1$, $m_2$ this is essentially the only choice. This was the result of extensive
symbolic computations carried out beginning at the time that~\cite{Gr} was written. This statement was proved in general in \cite{IT} a paper
that kindly acknowledges this work carried out initially in a special case.
If one allows the eight nearest
neighbours we get another linearly independent dif\/ference operator, a fact also found by symbolic computation by us and independently proved in \cite{IT}.

There is by now a rather large literature dealing with orthogonal polynomials in several variables. A reference that is still useful is \cite[Vol.~2]{HTF}. A comprehensive treatment is found in~\cite{DX}. For some of the issues that we are interested in the reader can see~\cite{IX,GI} and the references in these papers.

It may be useful to point out that when the conditions~\eqref{eq1.10} are not met the polynomials~$P_{m,n}$ still satisfy dif\/ference equations in the
indices $(m,n)$, of the type given in~\cite{IX,GI}. When the conditions~\eqref{eq1.10} are met these recursions become much simpler in that they involve a smaller number of neighbouring indices. Having a recursion relation involving the smallest number of neighbours of the index $(m,n)$ might be
important in certain numerical implementations of these recursions as well as in potential signal processing applications of these polynomials. In those situations, having a minimal number of sampling points could be a useful feature.

In dealing with the same polynomials Iliev and Terwilliger, see
the very recent nice paper~\cite{IT}, have found two $7$-term recurrence relations. In fact a referee has kindly pointed out that our $5$-term relation~\eqref{eq1.11} can be derived by using a combination of these.

It is clear that both in \cite{IT} and \cite{M1} the replacement of the variables $u_{i,j}$ by the consideration of
a matrix with entries $1-u_{i,j}$, properly augmented, played a very important role. This matrix is
considered in~\cite{IT} and an interpretation is given in terms of Lie algebras. In \cite{M1} the author proves that the orthogonality
of the columns of this matrix with respect to a weight built out of the~$\eta_i$
is the appropriate extension of~\eqref{eq1.10}.

The $F_1^{(2)}$ notation used in the $2$-variable case \eqref{eq1.8} and more generally in the $n$-variable case~\eqref{eq1.6} is a ref\/lection of the fact that these are generalizations of the standard Appell--Kamp\'e de F\'eriet function
\begin{gather*}
F_1(a;b,b';c;x,y) = \sum_i \sum_j \frac {(a)_{i+j}(b)_i(b')_j}{i!j!(c)_{i+j}} x^iy^j.
\end{gather*}

A very useful integral representation of this $F_1$ function is the double integral
\begin{gather*}
\frac {\G(c)}{\G(b)\G(b')\G(c-b-b')} \int_0^1 \!\int_0^1 \xi_1^{b-1}\xi_2^{b'-1}(1\!-\!\xi_1\!-\!\xi_2)^{c-b-b'-1}(1\!-\!\xi_1x\!-\!\xi_2y)^{-a}d\xi_1 d\xi_2,
\end{gather*}
provided
\[
0 < \mathrm{Re}(b,b',c-b-b').
\]
This extends to $F_1^{(2)}$ as well, which can be easily verif\/ied:
\begin{gather}
F_1^{(2)}(a_1,a_2;b_1,b_2;c;u_1,v_1,u_2,v_2)
 = \frac {\G(c)}{\G(a_1)\G(a_2)\G(c-a_1-a_2)}  \nonumber\\
  \qquad {}\x  \int_0^1 \int_0^1 \xi_1^{a_1-1}\xi_2^{a_2-1} (1-\xi_1-\xi_2)^{c-a_1-a_2-1}(1-u_1\xi_1 - u_2\xi_2)^{-b_1} \nonumber\\
\qquad {} \x (1-v_1\xi_1-v_2\xi_2)^{-b_2}d\xi_1d\xi_2,\label{eq1.14}
\end{gather}
which we shall f\/ind very useful in our calculations, even though the parameters in the case of~\eqref{eq1.8} do not satisfy the convergence conditions of the integral in \eqref{eq1.14}.  We will take the point of view that whatever identities we f\/ind by using~\eqref{eq1.14} with $0 < \mathrm{Re}(a_1,a_2,c-a_1-a_2)$, are also valid where $a_1$, $a_2$, $c$ are, in fact, negative integers.

For the direct hypergeometric proof that we are planning to give in the following pages it will be necessary to make use of the transformation formulas:
\begin{gather}
F_1^{(2)}(a_1,a_2;b_1,b_2;c;u_1,v_1,u_2,v_2) \label{eq1.15} \\
 \qquad {} = (1-v_1)^{-a_1}(1-v_2)^{-a_2}F_1^{(2)}\!\!\left(a_1,a_2;b_1,c\!-\!b_1\!-\!b_2;c;\frac {u_1\!-\!v_1}{1\!-\!v_1},\frac {-v_1}{1\!-\!v_1},\frac {u_2\!-\!v_2}{1\!-\!v_2},\frac {-v_2}{1\!-\!v_2}\right) \nonumber\\
\qquad {}= (1-u_1)^{-a_1}(1-u_2)^{-a_2} F_1^{(2)}\!\!\left(a_1,a_2;c\!-\!b_1\!-\!b_2,b_2;c;\frac {-u_1}{1\!-\!u_1},\frac {v_1\!-\!u_1}{1\!-\!u_1},\frac {-u_2}{1\!-\!u_2},\frac {v_2\!-\!u_2}{1\!-\!u_2}\right),\nonumber
\end{gather}
which were proved in Hoare and Rahman \cite{HR4}.  But there is a third transformation that we shall f\/ind occasions to use, that is valid when $(m_1,m_2)$ and $(x_1,x_2)$ are pairs of nonnegative integers, as is $N$, satisfying the triangle inequality: $0 \le m_1 + m_2 \le N$, $0 \le x_1 + x_2 \le N$, and that is
\begin{gather}
F_1^{(2)}(-m_1,-m_2;-x_1,-x_2;-N;u_1,v_1,u_2,v_2)   = \frac {(x_1+x_2-N)_{m_1+m_2}}{(-N)_{m_1+m_2}} \label{eq1.16}\\
{}\x F_1^{(2)}(-m_1,-m_2;-x_1,-x_2;N+1-x_1-x_2-m_1-m_2;1-u_1,1-v_1,1-u_2,1-v_2),\nonumber
\end{gather}
which is just a generalization of the transformation:
\begin{gather}
{}_2F_1(-m,-x;-N;u)  = \frac {(x-N)_m}{(-N)_m} {}_2F_1(-m,-x;N+1-x-m;1-u) \nonumber\\
\qquad = \frac {(m-N)_x}{(-N)_x} {}_2F_1(-m,-x;N+1-x-m;1-u).\label{eq1.17}
\end{gather}
In fact \eqref{eq1.16} and \eqref{eq1.17} easily extend to the multidimensional case $F_1^{(n)}$, provided one is dealing with terminating series. It may be remarked here that for \eqref{eq1.16} and \eqref{eq1.17} to be true, indeed in the general case of $F_1^{(n)}$, the parameter $N$ need not even be an integer.

\section{A general expression for orthogonality sum and proof of (\ref{eq1.10})}
\label{sec2}

Let us denote
\begin{gather*}
I_{m_1,m_2}^{n_1,n_2} = \underset{0 \le x_1 + x_2 \le N}{\underset{x_1,x_2}{\sum\sum}} b_2(x_1,x_2;N;\eta_1,\eta_2)P_{m_1,m_2}(x_1,x_2)P_{n_1,n_2}(x_1,x_2).
\end{gather*}
At $(m_1,m_2) = (0,0)$, and $(n_1,n_2) \ne (0,0)$, this simply represents a generating function for these polynomials, namely:
\begin{gather*}
I_{0,0}^{n_1,n_2} = (1-\eta_1u_1 - \eta_2v_1)^{n_1}(1-\eta_1u_2 - \eta_2v_2)^{n_2}.
\end{gather*}
So, at the three points $(0,0)$, $(1,0)$ and $(0,1)$ the pairwise orthogonality between the f\/irst and the last two simply amounts to the conditions (a) and (b) given in \eqref{eq1.10}.  To obtain the condition at the points $(1,0)$ and $(0,1)$ we need some more computations.

For ease of computation we will imagine, for the time being, that $-m_1$, $-m_2$, $-N$ are complex numbers $a_1$, $a_2$, $c$ such that $0 \le \mathrm{Re}(a_1,a_2,c-a_1-a_2)$.  Then using \eqref{eq1.14} we get
\begin{gather*}
I_{m_1,m_2}^{n_1,n_2}  = \frac {\G(c)}{\G(a_1)\G(a_2)\G(c-a_1-a_2)} \int_0^1 \int_0^1 d\xi_1d\xi_2\xi_1^{a_1-1}\xi_2^{a_2-1}(1-\xi_1-\xi_2)^{c-a_1-a_2-1} \\
\phantom{I_{m_1,m_2}^{n_1,n_2}  =}{}
\x\! \sum_{x_1}\sum_{x_2} b_2(x_1,x_2;N;\eta_1,\eta_2)P_{n_1,n_2}(x_1,x_2)
  (1-\xi_1u_1 - \xi_2u_2)^{x_1}(1-\xi_1v_1 - \xi_2v_2)^{x_2}.\nonumber
\end{gather*}
Clearly
\begin{gather}
\underset{x_1,x_2}{\sum\sum} b_2(x_1,x_2;N;\eta_1,\eta_2)(-x_1)_{i+k}(-x_2)_{j+l}  (1-\xi_1u_1 - \xi_2u_2)^{x_1}(1-\xi_1v_1 - \xi_2v_2)^{x_2}\nonumber \\
 \qquad{} = (-N)_{i+j+k+l}\eta_1^{i+k}\eta_2^{j+l}(1-\xi_1u_1 - \xi_2u_2)^{i+k}(1-\xi_1v_1 - \xi_2v_2)^{j+l} \nonumber\\
 \qquad \quad{} \x \{1 - \xi_1(\eta_1u_1 + \eta_2v_1) - \xi_2(\eta_1u_2 + \eta_2v_2)\}^{N-i-j-k-l} \nonumber\\
 \qquad {} = (-N)_{i+j+k+l}(\eta_1(1-\xi_1u_1 - \xi_2u_2))^{i+k}(\eta_2(1-\xi_1v_1 - \xi_2v_2))^{j+l},\label{eq2.4}
\end{gather}
by virtue of \eqref{eq1.10}(a) and \eqref{eq1.10}(b).  For general $(m_1,m_2)$ and $(n_1,n_2)$, we use \eqref{eq2.4} and recast back to the original parameters, getting
\begin{gather}
\label{eq2.5}
I_{m_1,m_2}^{n_1,n_2} = \sum \frac {(-n_1)_{i+j}(-n_2)_{k+l}}{i!j!k!l!} (\eta_1u_1)^i(\eta_2v_1)^j(\eta_1u_2)^k(\eta_2v_2)^l \\
\phantom{I_{m_1,m_2}^{n_1,n_2} =}{} \x \frac {(-i\!-\!j\!-\!k\!-\!l)_{m_1+m_2}}{(-N)_{m_1+m_2}}  F_1^{(2)} (-m_1,-m_2;-i\!-\!k,-j\!-\!l;-i\!-\!j\!-\!k\!-\!l;u_1,v_1,u_2,v_2).\nonumber
\end{gather}
Let us take $(m_1,m_2) = (0,1)$ and $(n_1,n_2) = (1,0)$, so that
\begin{gather*}
I_{0,1}^{1,0} = \sum_{i,j} \frac {(-1)_{i+j}}{i!j!} (\eta_1u_1)^i(\eta_2v_1)^j  \frac {(-i-j)}{(-N)} F_1(-1;-i,-j;-i-j;u_2,v_2) \nonumber\\
\phantom{I_{0,1}^{1,0}}{}
= \frac {(1-v_2)}{N} \sum_{i,j} \frac {(-1)_{i+j}}{i!j!} (\eta_1u_1)^i(\eta_2v_1)^j
 (i+j){}_2F_1\left[ \begin{matrix} -1,-i \\ -i-j \end{matrix}; \frac {u_2-v_2}{1-v_2}\right] \nonumber\\
\phantom{I_{0,1}^{1,0}}{}
= \frac {1-v_2}{N} \sum_{i,j} \frac {(-1)_{i+j}(\eta_1u_1)^i(\eta_2v_1)^j}{i!j!} j_2F_1\left[ \begin{matrix} -1,-i \\ j \end{matrix}; \frac {1-u_2}{1-v_2}\right] \nonumber\\
\phantom{I_{0,1}^{1,0}}{}
= \frac {1-v_2}{N} \sum_{i,j} \frac {(-1)_{i+j}}{i!j!} (\eta_1u_1)^i(\eta_2v_1)^j \left( j + i \frac {1-u_2}{1-v_2}\right) \nonumber\\
\phantom{I_{0,1}^{1,0}}{}
= \frac {1-v_2}{N} \left( \eta_2v_1 + \eta_1u_1 \frac {1-u_2}{1-v_2} \right)  = (\eta_1u_1(1-u_2) + \eta_2v_1(1-v_2))/N 
\end{gather*}
which must vanish, so \eqref{eq1.10}(c) must be satisf\/ied in addition to \eqref{eq1.10}(a) and \eqref{eq1.10}(b).

It is worth noting that by solving the f\/irst two conditions of \eqref{eq1.10} one can show that the third condition amounts to
\begin{gather}
\label{eq2.7}
U_1V_2 = U_2V_1,
\end{gather}
where
\[
U_i = 1 - u_i^{-1}, \qquad V_i = 1 - v_i^{-1}  , \qquad i = 1,2.
\]

Condition \eqref{eq2.7} has a simple geometrical interpretation as a cone embedded in four dimensional space. In a subsequent paper we will look
at a geometrical interpretation for the corresponding orthogonality conditions in the case of more than two variables.

We would also like to point out that if $\overline{\eta}_1$ and $\overline{\eta}_2$ are parameters of $b_2(m_1,m_2;N;\overline{\eta}_1,\overline{\eta}_2)$ for the dual orthogonality of the $P$'s, then they must satisfy
\begin{gather}
\text{(a)} \quad \overline{\eta}_1u_1 + \overline{\eta}_2u_2 = 1,\nonumber \\
\text{(b)}\quad \overline{\eta}_1v_1 + \overline{\eta}_2v_2 = 1, \label{eq2.8}\\
\text{(c)} \quad \overline{\eta}_1u_1v_1 + \overline{\eta}_2u_2v_2 = 1.\nonumber
\end{gather}

\section[Reduction of $I_{m_1,m_2}^{n_1,n_2}$]{Reduction of $\boldsymbol{I_{m_1,m_2}^{n_1,n_2}}$}
\label{sec3}

By the f\/irst transformation in \eqref{eq1.15} $F_1^{(2)}$ inside the sum in \eqref{eq2.5} becomes a multiple of $F_1$, which, transformed by \eqref{eq1.16} gives
\begin{gather*}
I_{m_1,m_2}^{n_1,n_2}  = (1-v_1)^{m_1}(1-v_2)^{m_2} \sum_{i,j,k,l} \frac {(-n_1)_{i+j}(-n_2)_{k+l}}{i!j!k!l!} (\eta_1u_1)^i(\eta_2v_1)^j(\eta_1u_2)^k(\eta_2v_2)^l \\
\phantom{I_{m_1,m_2}^{n_1,n_2}  =}{}
\x \frac {(-j-l)_{m_1+m_2}}{(-N)_{m_1+m+2}} F_1\left( -i-k;-m_1,-m_2;j+l+1-m_1-m_2;\frac {1-u_1}{1-v_1}, \frac {1-u_2}{1-v_2}\right).\nonumber
\end{gather*}
Set $i+k = r$, $j+l = s$, $i = r-k$, $j = s-l$, to get
\begin{gather}
I_{m_1,m_2}^{n_1,n_2} = (1-v_1)^{m_1}(1-v_2)^{m_2} \sum_{r,s} \frac {(-n_1)_{r+s}}{r!s!} (\eta_1u_1)^r(\eta_2v_1)^s \nonumber\\
\phantom{I_{m_1,m_2}^{n_1,n_2} =}{}
\x F_1\left( -n_2;-r,-s;n_1+1-r-s; \frac {u_2}{u_1}, \frac {v_2}{v_1}\right) \nonumber\\
\phantom{I_{m_1,m_2}^{n_1,n_2} =}{}
\x \frac {(-s)_{m_1+m_2}}{(-N)_{m_1+m_2}} F_1\left( -r;-m_1,-m_2;s+1-m_1-m_2; \frac {1-u_1}{1-v_1},\frac {1-u_2}{1-v_2}\right).\label{eq3.2}
\end{gather}
Since
\begin{gather*}
F_1\left(-n_2;-r,-s;n_1+1-r-s;\frac {u_2}{u_1},\frac {v_2}{v_1}\right) \nonumber\\
 \qquad {}= \frac {(-n_1-n_2)_{r+s}}{(-n_1)_{r+s}} F_1\left( -n_2;-r,-s;-n_1-n_2;1 - \frac {u_2}{u_1}, 1 - \frac {v_2}{v_1}\right)
\end{gather*}
by \eqref{eq1.16}, \eqref{eq3.2} reduces to
\begin{gather}
I_{m_1,m_2}^{n_1,n_2}  = \frac {(1-v_1)^{m_1}(1-v_2)^{m_2}}{(-N)_{m_1+m_2}} \sum_{r,s} \frac {(-n_1-n_2)_{r+s}}{r!s!} (\eta_1u_1)^r(\eta_2v_1)^s
\nonumber\\
 \phantom{I_{m_1,m_2}^{n_1,n_2}=}{} \x F_1\left( -n_2;-r,-s;-n_1-n_2;1 - \frac {u_2}{u_1}, 1 - \frac {v_2}{v_1}\right)
 \nonumber\\
 \phantom{I_{m_1,m_2}^{n_1,n_2}=}{} \x (-s)_{m_1+m_2}F_1 \left( -r;-m_1,-m_2;s+1-m_1-m_2; \frac {1-u_1}{1-v_1}, \frac {1-u_2}{1-v_2}\right).\label{eq3.4}
\end{gather}
To carry out the summations over $r$ and $s$ we employ the integral formula
\begin{gather*}
F_1(a;b,b';c;x,y)
 = \frac {\G(c)}{\G(a)\G(c-a)} \int_0^1 \xi^{a-1}(1-\xi)^{c-a-1}(1-\xi x)^{-b}(1-\xi y)^{-b'}d\xi,
\end{gather*}
see \cite[Vol.~1]{HTF}.  In our case $b = -r$, $b' = -s$, $x = 1 - \frac {u_2}{u_1}$, $y = 1 - \frac {v_2}{u_1}$, so in \eqref{eq3.4} we need to compute
\begin{gather}
\sum_{r,s} \frac {(-n_1-n_2)_{r+s}}{r!s!} \left( \eta_1u_1\left( 1 - \xi \left( 1 - \frac {u_2}{u_1}\right)\right)\right)^r \left( \eta_2v_1\left( 1 - \xi \left( 1 - \frac {v_2}{v_1}\right)\right)\right)^s (-r)_{i+j}\nonumber \\
 \qquad {} \x (-s)_{m_1+m_2-i-j}
 = (-1)^{m_1+m_2} (-n_1-n_2)_{m_1+m_2}\nonumber\\
 \qquad {} \x
 \left( \eta_1u_1 \left( 1 - \xi \left( 1 - \frac {u_2}{u_1}\right)\right)\right)^{i+j} \left( \eta_2v_1\left( 1 - \xi \left( 1 - \frac {v_2}{v_1}\right)\right)\right)^{m_1+m_2-i-j} \nonumber\\
\qquad{} \x \{1 - (1 - \xi)(\eta_1u_1 + \eta_2v_1) - \xi(\eta_1u_1u_2 + \eta_2v_1v_2)\}^{n_1+n_2-m_1-m_2},\label{eq3.6}
\end{gather}
with the implicit assumption that $n_1 + n_2 \ge m_1 + m_2$.  However, by \eqref{eq1.10}(a) and \eqref{eq1.10}(c) the expression in $\{\ \}$ vanishes unless $n_1 + n_2 = m_1 + m_2$.  Therefore, \eqref{eq3.6} becomes
\begin{gather*}
(m_1+m_2)!(\eta_2v_1)^{m_1+m_2}\left( \frac {\eta_1u_1}{\eta_2v_1}\right)^{i+j} \left( 1 - \xi \left( 1 - \frac {u_2}{u_1}\right)\right)^{i+j} \left( 1 - \xi \left( 1 - \frac {v_2}{v_1}\right)\right)^{m_1+m_2-i-j}\nonumber \\
\qquad {}\x \delta_{m_1+m_2,n_1+n_2},
\end{gather*}
and consequently,
\begin{gather}
I_{m_1,m_2}^{n_1,n_2} = \delta_{m_1+m_2,n_1+n_2}(1-v_1)^{m_1}(1-v_2)^{m_2}(\eta_2v_1)^{m_1+m_2} \frac {(m_1+m_2)!}{(-N)_{m_1+m_2}} \nonumber\\
\phantom{I_{m_1,m_2}^{n_1,n_2} =}{}
\x \sum_{i,j} \frac {(-m_1)_i(-m_2)_j}{i!j!} \left( -\frac {\eta_1u_1}{\eta_2v_1} \right)^{i+j} \left( \frac {1-u_1}{1-v_1} \right)^i \left( \frac {1-u_2}{1-v_2} \right)^j \nonumber\\
\phantom{I_{m_1,m_2}^{n_1,n_2} =}{}
\x F_1\left(-n_2;-i-j,i+j-m_1-m_2;-m_1-m_2; 1 - \frac {u_2}{u_1}, 1 - \frac {v_2}{v_1}\right) \nonumber\\
\phantom{I_{m_1,m_2}^{n_1,n_2}}{}
= \delta_{m_1+m_2,n_1+n_2}(1-v_1)^{m_1}(1-v_2)^{m_2}(\eta_2v_1)^{m_1+m_2} \left( \frac {u_2}{u_1} \right)^{n_2} \frac {(m_1+m_2)!}{(-N)_{m_1+m_2}} \nonumber\\
\phantom{I_{m_1,m_2}^{n_1,n_2} =}{}
\x \sum_{i,j} \frac {(-m_1)_i(-m_2)_j}{i!j!} \left( -\frac {\eta_1u_1}{\eta_2v_1} \right)^{i+j} \left( \frac {1-u_1}{1-v_1} \right)^i \left( \frac {1-v_2}{1-u_2} \right)^j \nonumber\\
\phantom{I_{m_1,m_2}^{n_1,n_2} =}{}
\x {}_2F_1\left( -n_2;i+j-m_1-m_2;-m_1-m_2;1 - \frac {u_1v_2}{u_2v_1} \right),\label{eq3.8}
\end{gather}
by a special case of the last identity \eqref{eq1.15}.

\section{Final summations in (\ref{eq3.8})}
\label{sec4}

At the last stage we will set $i+j = k$, $j = k-i$, so that the $i$-sum becomes
\begin{gather*}
{}_2F_1\left( -m_1,-k;m_2+1-k; \frac {(1-u_1)(1-v_2)}{(1-u_2)(1-v_1)}\right) \nonumber\\
 \qquad{} = \frac {(-m_1-m_2)_k}{(-m_2)_k} {}_2F_1 \left( -m_1,-k;-m_1-m_2; 1 - \frac {(1-u_1)(1-v_2)}{(1-u_2)(1-v_1)} \right)\nonumber \\
 \qquad {} = \frac {(-m_1-m_2)_k}{(-m_2)_k} {}_2F_1\left( -m_1,-k;-m_1-m_2; 1 - \frac {u_1v_2}{u_2v_1} \right),
\end{gather*}
by \eqref{eq2.7} and \eqref{eq1.17}.  Thus
\begin{gather*}
I_{m_1,m_2}^{n_1,n_2}  = \delta_{m_1+m_2,n_1+n_2} (1-v_1)^{m_1}(1-v_2)^{m_2}(\eta_2v_2)^{m_1+m_2} \frac {(m_1+m_2)!}{(-N)_{m_1+m_2}} \left( \frac {u_2}{u_1} \right)^{n_2}\\ 
 {}
 \x \sum_{k=0}^{m_1+m_2} \frac {(-m_1-m_2)_k}{k!} \sum_{i=0}^{m_1} \frac {(-m_1)_i(-k)_i}{i!(-m_1-m_2)_i} \sum_{j=0}^{n_2} \frac {(-n_2,k-m_1-m_2)_j}{j!(-m_1-m_2)_j}
  \left( 1 - \frac {u_1v_2}{u_2v_1} \right)^{i+j},\nonumber
\end{gather*}
since $-\frac {\eta_1u_1}{\eta_2v_1} \frac {(1-u_2)}{(1-v_2)} = 1$, because $\eta_1u_1(1-u_2) + \eta_2v_1(1-v_2) = 0$.  But, now
\begin{gather*}
\sum_{k=0} \frac {(-m_1-m_2)_k}{k!} (-k)_i(k-m_1-m_2)_j  \nonumber\\
\qquad {}  = (-1)^i (-m_1-m_2)_{i+j} \sum_{k=0}^{m_1+m_2-i-j} \frac {(i+j-m_1-m_2)_k}{k!} \nonumber\\
 \qquad {} = (-1)^{i} (-m_1-m_2)_{m_1+m_2} \delta_{m_1+m_2,i+j} \nonumber\\
 \qquad {} = (m_1+m_2)_{m_1+m_2} (-1)^{m_1+m_2-i} \delta_{m_1+m_2,i+j}.
\end{gather*}
So
\begin{gather*}
I_{m_1,m_2}^{n_1,n_2}  = \delta_{m_1+m_2,n_1+n_2} (1-v_1)^{m_1}(1-u_2)^{m_2} \left( \eta_2v_1 \left( 1 - \frac {u_1v_2}{u_2v_1} \right)\right)^{m_1+m_2} \left( \frac {u_2}{u_1} \right)^{n_2} \\
\phantom{I_{m_1,m_2}^{n_1,n_2}  =}{} \x \frac {((m_1+m_2)!)^2}{(-N)_{m_1+m_2}}   \sum_{i=0}^{m_1} \frac {(-m_1)_i(-1)^{m_1+m_2-i}}{i!(-m_1-m_2)_i} \frac {(-n_2)_{m_1+m_2-i}}{(-m_1-m_2)_{m_1+m_2-i}(m_1+m_2-i)!}.\nonumber
\end{gather*}
The summand is $0$ unless $n_2 \ge m_2 + m_1 - i \Rightarrow n_2 \ge m_2$, since $m_1 \ge i$.
So we set $i = m_1 + m_2 - n_2 + l$, $l \ge 0$, and get, for the $i$-sum above
\begin{gather*}
\frac {(-n_2)_{n_2}(-m_1)_{m_1+m_2-n_2}}{(-m_1-m_2)_{m_1+m_2}(m_1+m_2)!} \sum_{l=0}^{n_2-m_2} \frac {(m_2-n_2)_l}{l!} \nonumber\\
\qquad {} = \frac {(-m_2)_{m_2}(-m_1)_{m_1}}{(-m_1-m_2)_{m_1+m_2}(m_1+m_2)!} \delta_{m_2,n_2}
 = \frac {m_1!m_2!}{(m_1+m_2)!^2} \delta_{m_2,n_2} \quad \Rightarrow  \quad m_1 = n_1
\end{gather*}
since $m_1+m_2 = n_1+n_2$.  Thus,
\begin{gather*}
I_{m_1,m_2}^{n_1,n_2}  = \delta_{m_1,m_2}\delta_{n_1,n_2}(1-v_1)^{m_1}(1-v_2)^{m_2} \left( \frac {u_2}{u_1} \right)^{m_2} \left( -\eta_2v_1\left( 1 - \frac {u_1v_2}{u_2v_1} \right)\right)^{m_1+m_2} \nonumber\\
\phantom{I_{m_1,m_2}^{n_1,n_2}  =}{} \x 1\left/ \binom{N}{m_1,m_2}\right. .
 \end{gather*}
To determine the coef\/f\/icient in terms of $\overline{\eta}_1$ and $\overline{\eta}_2$, note that
\begin{gather*}
1 - \frac {u_1v_2}{u_2v_1}  = 1 - \frac {(1-u_1)(1-v_2)}{(1-u_2)(1-v_1)} = \frac {u_1-v_1+v_2-u_2 - (u_1v_2 - u_2v_1)}{(1-u_2)(1-v_1)} \nonumber\\
\hphantom{1 - \frac {u_1v_2}{u_2v_1}} {}= -\frac {D(1-\overline{\eta}_1-\overline{\eta}_2)}{(1-u_2)(1-v_1)},
\end{gather*}
from solving \eqref{eq2.8}, with
$D \equiv u_1v_2 - u_2v_1$.

Now, from \eqref{eq2.8}(a) and \eqref{eq2.8}(b) we get
\begin{gather*}
\overline{\eta_1} = (v_2-u_2)D^{-1},\qquad \overline{\eta_2} = (u_1-v_1)D^{-1},
\end{gather*}
while \eqref{eq1.10}(a)--\eqref{eq1.10}(c) give
\begin{gather*}
\eta_1 = \left. \begin{vmatrix} 1 & v_1 \\ 1 & v_1v_2 \end{vmatrix}\right/\begin{vmatrix} u_1 & v_1 \\ u_1u_2 & v_1v_2 \end{vmatrix}
= \frac {(v_2-1)}{u_1(v_2-u_2)} = -\frac {1-v_1}{u_2(v_1-u_1)}, \nonumber\\
\eta_2  = \frac {1-u_2}{v_1(v_2-u_2)} = \frac {1-u_1}{v_2(v_1-u_1)}.
\end{gather*}
Hence
\begin{gather*}
-\frac {D\eta_2v_1(1-v_1)}{(1-u_2)(1-v_1)} = -\frac {D}{v_2-u_2} = -\frac {1}{\overline{\eta}_1},
\end{gather*}
from \eqref{eq2.8}.

Now,
\[
\eta_2v_1(1-v_2) = -\eta_1u_1(1-u_2),
\qquad
\eta_2v_1(1-v_2) \frac {u_2}{u_1} = -\eta_1u_2(1-u_2),
\]
so
\begin{gather*}
-\frac {D\eta_2v_1(1-v_2)}{(1-u_2)(1-v_1)} \frac {u_2}{u_1} = \frac {D\eta_1u_2}{(1-u_1)} = -\frac {D}{v_1-u_1} = -\frac {1}{\overline{\eta}_2}.
\end{gather*}

Thus the normalization factor is
\[
\big(b_2(m_1,m_2;N;\overline{\eta}_1,\overline{\eta}_2)(1-\overline{\eta}_1-\overline{\eta}_2)^{-N}\big)^{-1}.
\]

\section{Proof of the recurrence relation (\ref{eq1.11})}
\label{sec5}

By the transformation formula \eqref{eq1.15},
\begin{gather*}
P_{m_1,m_2}(x_1,x_2)  = (1-u_1)^{x_1}(1-v_1)^{x_2}
  \sum \frac {(m_1+m_2-N)_{i+j}(-m_2)_{k+l}(-x_1)_{i+k}(-x_2)_{j+l}}{i!j!k!l!(-N)_{i+j+k+l}} \nonumber\\
\phantom{P_{m_1,m_2}(x_1,x_2)  =}{}
\x \left( \frac {-u_1}{1-u_1} \right)^i \left( \frac {-v_1}{1-v_1} \right)^j \left( \frac {u_2-u_1}{1-u_1} \right)^k \left( \frac {v_2-v_1}{1-v_1} \right)^l.
\end{gather*}
So, a straightforward calculation gives
\begin{gather*}
(N-m_1-m_2)(P_{m_1+1,m_2}(x_1,x_2) - P_{m_1,m_2}(x_1,x_2))
  = -(1-u_1)^{x_1}(1-v_1)^{x_2}\nonumber\\
  \qquad{}\x \left( u'_1 \frac {\p}{\p u'_1} + v'_1 \frac {\p}{\p v'_1} \right)((1-u_1)^{-x_1}(1-v_1)^{-x_2}P_{m_1,m_2}(x_1,x_2)),
\end{gather*}
where $u'_i = u_i/(u_i-1)$, $v'_i = v_i/(v_i-1)$, $i = 1,2$.

Clearly $u'_i \frac {\p}{\p u'_i} = u_i(1-u_i) \frac {\p}{\p u_i}$, etc.  Hence, with a similar expression for $P_{m_1,m_2+1} - P_{m_1,m_2}$, we can write
\begin{gather}
(N-m_1-m_2)\{A(P_{m_1+1,m_2} - P_{m_1,m_2}) - B(P_{m_1,m_2+1} - P_{m_1,m_2})\} \nonumber\\
 \qquad {} = B(1-u_2)^{x_1}(1-v_2)^{x_2} \left\{ u_2(1-u_2) \frac {\p}{\p u_2} + v_2 (1-v_2) \frac {\p}{\p v_2} \right\} \nonumber\\
 \qquad\quad {}  \x ((1-u_2)^{-x_1}(1-v_2)^{-x_2}P_{m_1,m_2})  - A(1-u_1)^{x_1}(1-v_1)^{x_2}\nonumber\\
\qquad\quad {}  \x
 \left\{ u_1(1-u_1) \frac {\p}{\p u_1} + v_1(1-v_1) \frac {\p}{\p v_1} \right\} ((1-u_1)^{-x_1}(1-v_1)^{-x_2}P_{m_1,m_2})
\label{eq5.3}
\end{gather}
for some suitably chosen constants $A$ and $B$.

A more convenient form of the right-hand side of \eqref{eq5.3} is
\begin{gather*}
\{x_1(Bu_2 - Au_1) + x_2(Bv_2 - Av_1)\}P_{m_1,m_2}
  + \left\{ B \left( \left(u_2 \left(1-u_2\right) \frac {\p}{\p u_2} + v_2\left( 1 - v_2\right) \frac {\p}{\p v_2}\right)\right)\right. \nonumber\\
 \qquad {} - A \left.\left( \left(u_1\left( 1 - u_1\right) \frac {\p}{\p u_1} + v_1\left( 1 - v_1\right) \frac {\p}{\p v_1} \right) \right) \right\} P_{m_1,m_2}.
\end{gather*}
By using the values of $A$ and $B$ indicated in \eqref{eq1.11}, and those of $u$'s and $v$'s in \eqref{eq1.4} we f\/ind that
\begin{gather*}
Bu_2 - Au_1 = p_1+p_2,\qquad Bv_2 - Av_1 = -(p_3+p_4).
\end{gather*}
Similarly,
\begin{gather*}
Cm_1(P_{m_1-1,m_2} - P_{m_1,m_2}) - Dm_2(P_{m_1,m_2-1}- P_{m_1,m_2}) \\
\qquad {}= \left\{ D \left( u_2 \frac {\p}{\p u_2} + v_2 \frac {\p}{\p v_2} \right) - C\left( u_1 \frac {\p}{\p u_1} + v_1 \frac {\p}{\p v_1} \right)\right\} P_{m_1,m_2}.
\end{gather*}
What we really need to prove is that
\begin{gather*}
\left\{ -(A(1-u_1)+C)u_1 \frac {\p}{\p u_1} - (A(1-v_1)+C)v_1 \frac {\p}{\p v_1} \right. \nonumber\\
\left.\qquad {} +  (B(1-u_2)+D)u_2 \frac {\p}{\p u_2} + (B(1-v_2)+D) v_2 \frac {\p}{\p v_2} \right\} P_{m_1,m_2}(x_1,x_2) = 0.
\end{gather*}
Straightforward algebra gives
\begin{alignat*}{3}
& A(1-u_1) + C  = p_4,\qquad && A(1-u_2) + C = -p_2, & \nonumber\\
& B(1-u_2) + D  = -p_3,\qquad && B(1-v_2) + D = p_1, & 
\end{alignat*}
so it amounts to showing that
\begin{gather}
\label{eq5.8}
\left( -p_4 u_1 \frac {\p}{\p u_1} + p_2v_1 \frac {\p}{\p v_1} - p_3u_2 \frac {\p}{\p u_2} + p_1v_2 \frac {\p}{\p v_2} \right) P_{m_1,m_2} = 0.
\end{gather}

Since the $u$'s and $v$'s are expressed in terms of the 4 $p$'s, what we need now is to express the derivatives in \eqref{eq5.8} in terms of those of the $p$'s.  Noting that
\begin{gather*}
(F_{u_1},F_{v_1},F_{u_2},F_{v_2})' = J^{-1}(F_{p_1},F_{p_2},F_{p_3},F_{p_4})',
\end{gather*}
for any dif\/ferentiable function $F$, with the Jacobian $J$ given by the $4 \x 4$ matrix
\begin{gather*}
J = (u_{1,j},v_{1,j},u_{2,j},v_{2,j}),
\qquad
(u_{i,j}) = \left( \frac {\p u_i}{\p p_1}, \frac {\p u_i}{\p p_2}, \frac {\p u_i}{\p p_3}, \frac {\p u_i}{\p p_4}\right)', \qquad \text{etc.}, \qquad i = 1,2,
\end{gather*}
we are reduced to the task of proving that
\begin{gather}
\label{eq5.12}
-p_4 u_1a_{11} + p_2v_1a_{21} - p_3u_2a_{31} + p_1v_2a_{41} = 0,
\end{gather}
and 3 more similar relations, where the $|J|^{-1}a'_{ij}$ are elements of the inverse matrix $J^{-1}$, which, of course, exists.  By a set of long and messy computations we obtain
\begin{alignat}{3}
& a_{11} = -p_1v_2\D^2/p_4p_2^2p_3^2S^3,\qquad && a_{21} = -u_2\D^2/p_1p_2p_4^2S^3, & \nonumber\\
& a_{31} = -v_1\D^2/p_1p_3p_4^2S^3,\qquad &&  a_{41} = -u_1\D^2/p_2^2p_3^2S^3,& \label{eq5.13}
\end{alignat}
with $S = p_1+p_2+p_3+p_4$, $\D = p_1p_4 - p_2p_3$.  Substitution of \eqref{eq5.13} proves \eqref{eq5.12}.  The three other relations are similarly proved.
This completes the proof of \eqref{eq1.11}.

\medskip

{\bf An additional comment.} After this paper was completed we became aware of a recent arXiv posting~\cite{IT},
where the authors point out some important work of H.~Mizukawa and H.~Tanaka~\cite{MT}.
In a future publication we return to the probabilistic origin of the work of M.~Hoare and M.~Rahman
and we discuss the relation between the approach in~\cite{MT}, based on the notion of character
algebras, and our own.

\subsection*{Acknowledgements}

We thank one referee in particular for a very methodical job that has rendered this into a more
accurate paper. In an area where several people have made important contributions he has helped us
tell the story properly.

The research of the f\/irst author was supported in part by the Applied Math. Sciences subprogram of the Of\/f\/ice of Energy Research, USDOE, under Contract DE-AC03-76SF00098, and by AFOSR under contract FA9550-08-1-0169.

\pdfbookmark[1]{References}{ref}
\LastPageEnding

\end{document}